\documentclass[11pt]{amsart}
\usepackage{amsmath}
\usepackage{amsfonts}
\usepackage{amssymb}
\usepackage{enumerate}
\usepackage{amsthm}
\usepackage{graphicx}
\usepackage{verbatim}

\theoremstyle{definition}
\newtheorem{proposition}{Proposition}[section]
\newtheorem{theorem}[proposition]{Theorem}

\newtheorem{corollary}[proposition]{Corollary}

\newtheorem{example}[proposition]{Example}
\newtheorem{remark}[proposition]{Remark}

\begin{document}
\title[Error bounds for multivariate Laplace approximation]{Simple error bounds for the multivariate Laplace approximation under weak local assumptions}

\author{P. Majerski}
\address{Piotr Majerski, Faculty of Applied Mathematics,
						 AGH University of Scien\-ce and Technology,
						 Al. Mickiewicza 30, 30-059 Krak\'ow, Poland}
\email{majerski@agh.edu.pl}

\maketitle

\begin{abstract}
The paper provides new upper and lower bounds for the multivariate Laplace approximation under weak local assumptions. Their range of validity is also given. An application to an integral arising in the extension of the Dixon's identity is presented. The paper both generalizes and complements recent results by Inglot and Majerski and removes their superfluous assumption on vanishing of the third order partial derivatives of the exponent function. 
\end{abstract}

\section{Notation and preliminaries}

%

For a given function $f:\mathbb{R}^d\to\mathbb{R}$, with $d\geqslant2$, $\dot{f}(\mathbf{0})$ and 
$\ddot{f}(\mathbf{0})$ is used to denote, respectively, the gradient vector and the Hessian matrix of the function $f$ evaluated 
at $\mathbf{0}=(0,\ldots,0)^\prime\in\mathbb{R}^d$. $d^2f(\boldsymbol0,\boldsymbol t)$, $d^3f(\boldsymbol0,\boldsymbol t)$ and $d^4f(\boldsymbol0,\boldsymbol t)$ will respectively denote the second, the third and the fourth order total derivative of $f$ at $\mathbf{0}$ with respect to $\mathbf{t}=(t_1,\ldots,t_d)^\prime\in\mathbb{R}^d$, so in particular $d^2f(\boldsymbol0,\boldsymbol t)=\boldsymbol t'\ddot{f}(\mathbf{0})\boldsymbol t$, and
\[
	d^3f(\mathbf{0},\mathbf{t}):=\sum_{i=1}^d\sum_{j=1}^d\sum_{k=1}^d \frac{\partial^3f(\mathbf{0})}{\partial t_i\partial t_j\partial t_k}t_it_jt_k.
\]
For $r>0$ and a positive definite $d\times d$ matrix $A$, we define $B_r:=\{\mathbf{t}\in\mathbb{R}^d:~\mathbf{t}^\prime\mathbf{t}\leqslant r^2\}$ and $E_r(A):=\{\mathbf{t}\in\mathbb{R}^d:~\mathbf{t}^\prime A\mathbf{t}\leqslant r^2\}$.

Finally, concerning the smoothness of functions, for a nonnegative integer $k$ by $f\in {\rm C}^k(\Omega)$ or $f\in {\rm C}^k(\boldsymbol t_0)$ we mean that $f$ is of class ${\rm C}^k$ on $\Omega$ or some neighborhood of the point $\boldsymbol t_0$, respectively. We shall use an analogous notation for the H\"older continuous class of functions ${\rm C}^{k,\alpha}$, with $\alpha\in(0,1]$.

For $n>0$ and an integer $d\ge2$ consider the Laplace integral
\[
	J(n):=\int_{\Omega}e^{-nf(\boldsymbol t)}g(\boldsymbol t)d\boldsymbol t,
\]
where $\Omega\subset\mathbb{R}^d$ is a measurable set and $f,g$ are two real valued measurable functions defined on $\Omega$. Under classical assumptions, i.e.
\begin{itemize}
	\item [(L1)] $J(n)$ exists and is finite for every large $n$,
	\item [(L2)] $f$ attains a separate absolute minimum value at an interior point $\boldsymbol t_0$ of $\Omega$,
	\item [(L3)] $f\in {\rm C}^2(\boldsymbol t_0)$, $\dot{f}(\boldsymbol t_0)=\boldsymbol 0$, with $\boldsymbol0=(0,\ldots,0)'\in\mathbb{R}^d$; $\ddot{f}(\boldsymbol t_0)>0$,
	\item [(L4)] $g$ is continuous at $t_0$ and $g(t_0)\neq0$,
\end{itemize}
it is well known (\cite{BH,deB,H,Hsu1,W}) that as $n\to\infty$, $J(n)$ is asymptotic to $\tilde{J}(n)$, where
\[
	\tilde{J}(n):=\frac{e^{-n f(\boldsymbol t_0)}g(\boldsymbol t_0)}{\sqrt{\det(\ddot{f}(\boldsymbol t_0))}}\left(\frac{2\pi}{n}\right)^\frac{d}{2}.
\]
There are two broad issues concerning the quality of this approximation, namely the rate at which its relative error 
	$$E(n):=J(n)/\tilde{J}(n)-1$$
converges to 0 (inseparably with the associated regularity of $f$ and $g$) and the expression for the coefficient which appears in the error term. Let us first turn to the second question and, for the moment, notice only that typically this rate is $n^{-1}$. Several authors devoted their efforts to evaluating this constant (as well as these appearing in a further asymptotic expansion), always encountering an inevitable complexity in the obtained formulae, no matter if they considered a general setup as Kirwin (\cite{K}) or a special form of the exponent function $f$ as Kaminski and Paris (\cite{KP1,KP2}). This is actually unavoidable even in the one dimensional case (see, e.g., \cite{deB,Wo1,Wo2,LPPS}), which is not considered in this paper.

From a practical point of view, rather than an asymptotic constant, one needs the upper and lower bounds for $E(n)$ instead, valid for $n\ge n_0$, say, with a benefit of a knowledge of $n_0$. The upper bounds for $|E(n)|$ was first derived in the paper \cite{MCW}, where the two dimensional case ($d=2$) was considered in detail. Under global regularity conditions ($f\in {\rm C}^4(\Omega)$, $g\in{\rm C}^2(\Omega))$) and assuming some restrictions on the rate of growth of $f$ over $\Omega$ McClure and Wong provided such bounds valid for every $n$ for which $J(n)$ exists, and proved the typical rate $n^{-1}$ of the relative error of Laplace approximation. The main advantage of this result is its range of validity. On the other hand, even for $d=2$, it gives complicated coefficients of the main error term.

Recently, in \cite{IM} Inglot and Majerski introduced a new proof the Laplace approximation and obtained for any $d\ge2$ simple upper and lower bounds for $E(n)$, which essentially depend the smallest eigenvalue of the Hessian matrix $\ddot{f}(\boldsymbol t_0)$ only.
Despite this result requires a weak local regularity, $f\in{\rm C}^{3,1}(\boldsymbol t_0)$, it holds only for a constant $g$ and
under the hypotheses that the third order partial derivatives of $f$ all vanish at $\boldsymbol t_0$. Especially the last condition is constricting and was conjectured superfluous.

In the presented paper we shall prove that similar results actually hold without these constraints, i.e. when just $f\in{\rm C}^{3,1}(\boldsymbol t_0)$ and $g\in{\rm C}^{1,1}(\boldsymbol t_0)$. We also pay attention to make our results easier to use and remove two technical parameters used in \cite{IM}. Thus, up to these minor convenience aimed changes the present results do generalize their counterparts from \cite{IM}.

To complete the discussion on the rate of convergence of $E(n)$ it is worth mentioning here, that the typical rate $n^{-1}$ is not the only possible and is slower for less regular functions $f$ and $g$, but improvement of their regularity is insufficient alone to speed it up. In fact it follows from \cite{IM} that:
\begin{itemize} 
	\item if $f\in {\rm C}^{2,\alpha}(\boldsymbol t_0)$, $g\in{\rm C}^{0,\beta}(\boldsymbol t_0)$ for some $\alpha,\beta\in(0,1]$ then $E(n)$ of the order $n^{-\alpha\wedge \beta}$;
	\item if $f\in{\rm C}^{k,\alpha}(\boldsymbol t_0)$ with $k+\alpha>4$, $g$ is a constant function and when for $i=3,\ldots,k-1$ all $i$th order partial derivatives vanish at $\boldsymbol t_0$, then $E(n)$ is of the order $n^{-1-(k+\alpha-4)/2}$.
\end{itemize}

An application of our results to previously inaccessible examples considered earlier in \cite{H} and \cite{MCW} is presented. In particular we find that for $d=2$ our estimations give a better main error term than those obtained by the McClure and Wong technique and it turns out that a range of a factual validity of our bounds is far better than that indicated by our theoretical result. The latter is inevitable in the Inglot and Majerski approach and fortiori in the present one. We, however discuss a remedy to this flaw by sacrificing the quality of the constant in the main error term.

Finally, in addition to the presented application we point out that a motivation to study the foundations of the Laplace expansion comes from different fields of mathematical statistics (for some recent papers see \cite{NFL,MS,MS2,EZS,BC,I,MR}) but is also useful in such branches of science as natural language processing (see, e.g., \cite{HA-MJ} and the references therein) and statistical mechanics and quantum field theory (cf, \cite[p. 232]{K}).

\section{Main result}
We first deal with integrals of a fully exponential form, i.e. 
$$I(n):=\int_{\Omega}e^{-n f(\mathbf{t})}d\mathbf{t}.$$ 
For simplicity we shall consider the normalized situation, which means that $f$ attains its global minimum $\mathbf{t_0}$ at $\mathbf{0}\in\mathbb{R}^d$ and $f(\mathbf{0})=0$. This is, of course, no restriction since putting $f_0(\mathbf{t}):=f(\mathbf{t}+\mathbf{t_0})-f(\mathbf{t_0})$ and $\Omega_0:=\{\mathbf{t}-\mathbf{t_0},~~{\mathbf t}\in \Omega\}$ one can write  
\[
	\int_{\Omega}e^{-n f(\mathbf{t})}d\mathbf{t}=e^{-n f(\mathbf{t_0})}\int_{\Omega_0}e^{-n f_0(\mathbf{t})}d\mathbf{t}.
\]

In what follows $||\cdot||$ will denote the Euclidean norm.

\begin{theorem}\label{th:muldimlapl}
    Let $\Omega\subset \mathbb{R}^d,\;d\geqslant 2$ be a measurable set, let $f:\Omega\to\mathbb{R}$ be a measurable function, and suppose that for some $n_1>0$ the integral
    $I(n)=\int_{\Omega}e^{-n f(\mathbf{t})}d\mathbf{t}$ exists and is finite for all $n\geqslant n_1$. Let $\dot f$ vanish at $\mathbf{0}$, $\ddot f(\mathbf{0})$ exist and be positive definite with $\lambda_{\rm{min}}$ being its smallest eigenvalue, all third order partial derivatives of $f$ exist at $\boldsymbol0$. Denote
    $$
	D:=\frac{d^{3/2}}6\max_{i,j,k}\left|\frac{\partial ^3f(\mathbf{\mathbf{0}})}{\partial t_i\partial t_j\partial t_k}\right|,\qquad i,j,k\in\{1,\ldots,d\}.
	$$ 
  Assume that\\

    {\rm (A1)}\quad there exist positive real numbers $r$, $C$ and $\alpha>1$ such that $B_r\subset\Omega$ and for every $\mathbf{t}\in B_r$
    $$\left|f(\mathbf{t})-\frac12d^2f(\mathbf{0},\mathbf{t})-\frac16d^3f(\mathbf{0},\mathbf{t})\right|\leqslant C ||\mathbf{t}||^{2+\alpha};
    $$

    {\rm (A2)}\quad there exist $\delta>0$ and $\Delta>0$ such that $f$ is convex on $B_{\delta}\subset\Omega$ and, moreover, for every $\mathbf{t}\in \Omega\setminus B_{\delta\wedge r}$
    $$f(\mathbf{t})\geqslant \Delta.
    $$
	    
	Then for every $n\geqslant n_0$ with $n_0$ given explicitly below, the following bounds hold true
    \begin{equation}\label{Lap_appr_statement_Lower}
    I( n)\geqslant\frac{1}{\sqrt{\mbox{\rm det} \ddot{f}(\mathbf{0})}}
    \left(\frac{2\pi}{ n}\right)^{d/2}\left[1-\frac{K_{\alpha,1}}{ n^{\alpha/2}}-\frac{K_l}{ n^{1+\alpha}}\right],
    \end{equation}
    \begin{equation}\label{Lap_appr_statement_Upper}
    I( n)\leqslant\frac{1}{\sqrt{\mbox{\rm det} \ddot{f}(\mathbf{0})}}
    \left(\frac{2\pi}{ n}\right)^{d/2}\left[1+\frac{K_{\alpha,1}}{ n^{\alpha/2}}+\frac{K_1}{n}+\frac{K_{\alpha,2}+K_u}{ n^{\alpha}}\right],
    \end{equation}
    where
    \begin{equation}\label{constant_K}
        K_{\alpha,1}=C(2/\lambda_{\rm{min}})^{1+\alpha/2}(d/2)_{1+\alpha/2},\quad 
		K_{\alpha,2}=C^2(2/\lambda_{\rm{min}})^{2+\alpha}(d/2)_{2+\alpha},
    \end{equation}
    \begin{equation}\label{constant_K_1}
	   K_1=D^2(2/\lambda_{\rm{min}})^3(d/2)_3,
    \end{equation}
    \begin{equation}\label{constant_K_l}
        K_l=\frac{e}2\sqrt{\frac{d}\pi}\left(1+\frac{2\alpha}d\right)^{d/2-1},
    \end{equation}
        \begin{equation}\label{constant_K_u}
		K_u=\frac{7\sqrt{\mbox{\rm det}\ddot{f}(\mathbf{0})}}{4(2\pi)^{d/2}}I( n_1)e^{\xi  n_1/2},
    \end{equation}
    $\xi=(r^2\lambda_{\rm{min}})\wedge (2\Delta)$ and $(x)_a=\Gamma(x+a)/\Gamma(x)$ is the Pochhammer symbol ($a\ge0$).
    
    Moreover, $n_0:=\inf (N_1\cap  N_2)$, with
    \[
		N_1:=\left\{
				n\ge1:~{d}\leqslant(d+2\alpha){\log n}\leqslant \xi n 
			 \right\},	
    \]
    \[
		N_2:=\left\{
				n\ge1:~D\left(\frac{d+2\alpha}{\lambda_{\rm{min}}}\right)^{3/2}\cdot\frac{\log^{3/2}n}{n^{1/2}}+ C\left(\frac{d+2\alpha}{\lambda_{\rm{min}}}\right)^{1+\alpha/2}\cdot\frac{\log^{1+\alpha/2}n}{n^{\alpha/2}}\leqslant7/4
			 \right\}.	
    \]
\end{theorem}

\noindent{\bf Proof.} 
Without loss of generality one can set $\Omega:=\mathbb{R}^d$. In fact, the suitable estimations in the proof below are eventually done by integrating over $\mathbb{R}^d$. 

Let $U$ be the upper triangular matrix with positive diagonal entries defined as $\ddot{f}(\mathbf{0})=U^\prime U$. For $ n\geqslant 1$ put $\varepsilon=\sqrt{\frac{(d+2\alpha)\log n}{n}}$ and note that $\varepsilon=\varepsilon( n)\to0$ as $ n\to\infty$. Then $I( n)$ can be written as
\[
    I( n)=\int_{E_\varepsilon(\ddot{f}(\mathbf{0}))}e^{- n f(\mathbf{t})}d\mathbf{t}+\int_{(E_\varepsilon(\ddot{f}(\mathbf{0})))^c}e^{- n f(\mathbf{t})}d\mathbf{t}=I_1( n)+I_2( n),
\]
where ${(E_\varepsilon(\ddot{f}(\mathbf{0})))^c}:=\mathbb{R}^d\setminus {E_\varepsilon(\ddot{f}(\mathbf{0}))}$.
Let $ n_2\geqslant 1\vee n_1$ be the smallest number such that
\begin{equation}\label{eq:lambda2}
	\frac{d}{n}\leqslant(d+2\alpha)\frac{\log n}{ n}\leqslant \xi 
\end{equation}
for all $ n\geqslant  n_2$. It then follows from Rayleigh-Ritz theorem (\cite{HJ}) that for\\
$ n\geqslant n_2$ it holds $\varepsilon^2\leqslant \lambda_{\rm{min}} r^2$ and $E_\varepsilon(\ddot f(\mathbf{0}))\subset B_r$. Using (A1) we have
\begin{align}
    I_1( n)&=\int_{E_\varepsilon(\ddot{f}(\mathbf{0}))}e^{- n f(\mathbf{t})}d\mathbf{t}\notag\\
                &\leqslant \int_{E_\varepsilon(\ddot{f}(\mathbf{0}))}e^{- n\mathbf{t}^\prime\ddot{f}(\mathbf{0})\mathbf{t}/2 -nd^3f(\mathbf{0},\mathbf{t})/6 + C n ||\mathbf{t}||^{2+\alpha}}d\mathbf{t}\notag\\
                &=\frac{n^{-d/2}}{\det U} \int_{B_{\sqrt{ n}\varepsilon}}e^{-||\mathbf{u}||^2/2} \exp\left\{-\frac{d^3f(\mathbf{0},U^{-1}\mathbf{u})}{6n^{1/2}}+ \frac{C||U^{-1}\mathbf{u}||^{2+\alpha}}{ n^{\alpha/2}}\right\}d\mathbf{u}\label{ineq:I1U1},
\end{align}
where in (\ref{ineq:I1U1}) we have applied the substitution $\mathbf{u}=\sqrt{ n}U\mathbf{t}$.
The condition
\begin{equation}\label{eq:lambda0}
	D\left(\frac{d+2\alpha}{\lambda_{\rm{min}}}\right)^{3/2}\cdot\frac{\log^{3/2}n}{n^{1/2}}+ C\left(\frac{d+2\alpha}{\lambda_{\rm{min}}}\right)^{1+\alpha/2}\cdot\frac{\log^{1+\alpha/2}n}{n^{\alpha/2}}\leqslant7/4
\end{equation}
insures that the quantity in curly brackets in (\ref{ineq:I1U1}) is less or equal $7/4$ for all $\mathbf{u}\in B_{\sqrt{n}\varepsilon}$. Let thus $n_0\geqslant n_2$ be the smallest number such that the inequality in (\ref{eq:lambda0}) holds for all $ n\geqslant  n_0$. Using the the relation $\|U^{-1}\mathbf{u}\|\leqslant \lambda_{\rm{min}}^{-1/2} \|\mathbf{u}\|$, the definition of $\varepsilon$ and the inequality $e^x\leqslant1+x+x^2$ ($x\in(-\infty,7/4]$), we get from (\ref{ineq:I1U1}) for $ n\geqslant  n_0$
\begin{align}
	I_1(n)&\leqslant\frac{n^{-d/2}}{{\rm det}\,U}\int_{B_{\sqrt{n}\varepsilon}}e^{-\|\mathbf{u}\|^2/2}
			\Bigg\{1+\frac{C\|\mathbf{u}\|^{2+\alpha}}{\lambda_{\rm{min}}^{1+\alpha/2}n^{\alpha/2}}-
					\frac{d^3f(\mathbf{0},U^{-1}\mathbf{u})}{6n^{1/2}}\notag\\
		  &+\frac{C^2\|\mathbf{u}\|^{4+2\alpha}}{\lambda_{\rm{min}}^{2+\alpha}n^\alpha}-
			\frac{C\|\mathbf{u}\|^{2+\alpha}d^3f(\mathbf{0},U^{-1}\mathbf{u})}{3n^{(1+\alpha)/2}\lambda_{\rm{min}}^{1+\alpha/2}}+
			\frac{\left(d^3f(\mathbf{0},U^{-1}\mathbf{u})\right)^2}{36n}\Bigg\}d\mathbf{u}.\notag
\end{align}
Observe that for every $R>0$, $\beta\geqslant0$ and for every $d\times d$ matrix $A$ we have
\begin{equation}\label{eq:int_d3}
	\int_{B_R}e^{-\|\mathbf{u}\|^2/2}\|\mathbf{u}\|^\beta d^3f(\mathbf{0},A\mathbf{u})d\mathbf{u}=0.
\end{equation}
Hence, using the inequality
\[
	\int_{B_{\sqrt{n}\varepsilon}}e^{-\|\mathbf{u}\|^2/2}d\mathbf{u}\leqslant (2\pi)^{d/2}
\]
we get for $n\geqslant n_0$
\begin{align}
\det U n^{d/2} I_1( n)&\leqslant (2\pi)^{d/2}+
\frac{C}{\lambda_{\rm{min}}^{1+\alpha/2} n^{\alpha/2}}\int_{\mathbb{R}^d} ||\mathbf{u}||^{2+\alpha}e^{-||\mathbf{u}||^2/2}d\mathbf{u}\notag\\
& +\frac{C^2}{\lambda_{\rm{min}}^{2+\alpha} n^{\alpha}}\int_{\mathbb{R}^d} ||\mathbf{u}||^{4+2\alpha}e^{-||\mathbf{u}||^2/2}d\mathbf{u}+\frac{D^2}{\lambda_{\rm{min}}^3n}\int_{\mathbb{R}^d} ||\mathbf{u}||^6e^{-||\mathbf{u}||^2/2}d\mathbf{u}\notag\\
&=(2\pi)^{d/2}\left[1+
\frac{K_{\alpha,1}}{ n^{\alpha/2}}+\frac{K_{\alpha,2}}{ n^{\alpha}}+\frac{K_1}{n}\right].\label{ineq:I1U}
\end{align}

Arguing in the same way, using again (A1), the inequality $e^x\geqslant 1+x$ and (\ref{eq:int_d3}) we obtain for $ n\geqslant n_0$
\begin{equation}\label{ineq:I1L}
    {\det U} n^{d/2}I_1( n)\geqslant \int_{B_{\sqrt{ n}\varepsilon}} e^{-||\mathbf{u}||^2/2}d\mathbf{u}-\left(2\pi\right)^{d/2}\frac{K_{\alpha,1}}{ n^{\alpha/2}}.
\end{equation}

Set $R^2=n\varepsilon^2$. Using Lemma 2.1 from \cite{IM} we have for $n\ge n_2$
\begin{align}
	\int_{B_{\sqrt{ n}\varepsilon}} e^{-||\mathbf{u}||^2/2}d\mathbf{u}&\ge (2\pi)^{d/2}\left(1-\frac{1}{\sqrt{\pi d}}\frac{R^2}{R^2-d+2}\exp\left\{-\frac12\left(R^2-d-(d-2)\log\frac{R^2}{d}\right)\right\}\right)\notag\\
	&=(2\pi)^{d/2}\left(1-n^{-(\alpha+1)}\frac{e^{{d}/2}}{\sqrt{\pi d}}\frac{R^2}{R^2-d+2}\exp\left\{-\frac{d-2}{2(d+2\alpha )}R^2+\left(\frac{d}2-1\right)\log\frac{R^2}{d}\right\}\right),\notag
\end{align}
where the equality $R^2=2(\alpha+1)\log n+(d-2)/(d+2\alpha)R^2$ has been used. Observing that $R^2/(R^2-d+2)$ decreases with $R^2$ while the exponent attains its maximum value at $R^2=d+2\alpha$ we get for $n\ge n_2$ the following lower bound for the integral $\int_{B_{\sqrt{ n}\varepsilon}} e^{-||\mathbf{u}||^2/2}d\mathbf{u}$,
\[
	(2\pi)^{d/2}\left(1-n^{-(\alpha+1 )}(e/2)\sqrt{d/\pi}(1+2\alpha/d)^{d/2-1}\right).
\]
This together with (\ref{ineq:I1L}) and the fact that $I_2( n)\geqslant 0$ proves the lower bound in (\ref{Lap_appr_statement_Lower}) for $ n\geqslant n_0$.

In order to prove the upper bound in (\ref{Lap_appr_statement_Upper}) one has to upper bound $I_2( n)$.
To this end observe first that from $({\rm A2})$ and (\ref{eq:lambda2}) we have on the set $(E_{\varepsilon}(\ddot{f}(\mathbf{0})))^c\cap B_{\delta\wedge r}^c$, for every $n\ge n_2$
\begin{equation}\label{eq:EcBc}
	e^{-(n-n_1)f(\mathbf{t})}\le e^{-n\Delta}e^{n_1\Delta}
							 \le n^{-\frac{d}2-\alpha}e^{n_1\xi/2}.
\end{equation}

Now let us turn to the set $(E_{\varepsilon}(\ddot{f}(\mathbf{0})))^c\cap B_{\delta\wedge r}$. For every $\mathbf{t}$ from this intersection define $\mathbf{t}^*=(\varepsilon/(\mathbf{t}'\ddot{f}(\mathbf{0})\mathbf{t})^{1/2})\cdot \mathbf{t}$ and note that 
$\mathbf{t}^*$ belongs to the line segment joining $\mathbf{0}$ with $\mathbf{t}$. Observe also that $(\mathbf{t}^*)'\ddot{f}(\mathbf{0})\mathbf{t}^*=\epsilon^2$.
 Convexity of $f$ together with the fact that $\mathbf{0}$ is its minimizer imply that $f(\mathbf{t})\ge f(\mathbf{t}^*)$, 
hence from the assumption (${\rm A1}$) for every $n\ge n_2$
\[
	f(\mathbf{t})\ge f(\mathbf{t}^*)\ge\frac12(\mathbf{t}^*)'\ddot{f}(\mathbf{0})\mathbf{t}^*+\frac16d^3f(\mathbf{0},\mathbf{t}^*)-C\|\mathbf{t}^*\|^{2+\alpha}
\]
on the set $(E_{\varepsilon}(\ddot{f}(\mathbf{0})))^c\cap B_{\delta\wedge r}$.
Since
\[
	\|\mathbf{t}^*\|\le\lambda_{{\rm min}}^{-1}(\mathbf{t}^*)'\ddot{f}(\mathbf{0})\mathbf{t}^*
\]
and
\[
	\left|\frac16d^3f(\mathbf{0},\mathbf{t}^*)\right|\le D\|\mathbf{t}^*\|^3\le
	\frac{D}{\lambda_{{\rm min}}^{3/2}}((\mathbf{t}^*)'\ddot{f}(\mathbf{0})\mathbf{t}^*)^{3/2}
\]
we get for every $n\ge n_2$
\begin{equation}\label{eq:fLowerEcB}
	f(t)\ge \frac12\varepsilon^2-\frac{D}{\lambda_{{\rm min}}^{3/2}}\varepsilon^3-\frac{C}{\lambda_{{\rm min}}^{1+\alpha/2}}\varepsilon^{2+\alpha}
\end{equation}
on the set $(E_{\varepsilon}(\ddot{f}(\mathbf{0})))^c\cap B_{\delta\wedge r}$. We also have for $n\ge n_0\ge n_2$
\begin{align}
	&\exp\left\{-(n-n_1)\left(\frac12\varepsilon^2-\frac{D}{\lambda_{{\mathrm min}}^{3/2}}\varepsilon^3 -\frac{C}{\lambda_{{\mathrm min}}^{1+\alpha/2}}\varepsilon^{2+\alpha}\right)\right\}\notag\\
	&=\exp\left\{\frac{n}2\varepsilon^2\right\}\exp\left\{n\left(\frac{D}{\lambda_{{\mathrm min}}^{3/2}}\varepsilon^3 +\frac{C}{\lambda_{{\mathrm min}}^{1+\alpha/2}}\varepsilon^{2+\alpha}\right)\right\}\exp\left\{n_1\left(\frac{\varepsilon^2}2-\frac{D}{\lambda_{{\mathrm min}}^{3/2}}\varepsilon^3 -\frac{C}{\lambda_{{\mathrm min}}^{1+\alpha/2}}\varepsilon^{2+\alpha}\right)\right\}\notag\\
	&\le \frac74n^{-({d}/2+\alpha)}\exp\left\{n_1\left(\frac{\varepsilon^2}2-\frac{D}{\lambda_{{\mathrm min}}^{3/2}}\varepsilon^3 -\frac{C}{\lambda_{{\mathrm min}}^{1+\alpha/2}}\varepsilon^{2+\alpha}\right)\right\}\notag\\
	&\le \frac74n^{-({d}/2+\alpha)}\exp\left\{n_1\frac{\varepsilon^2}2\right\}\le\notag\frac74e^{\xi n_1/2}n^{-({d}/2+\alpha)},
\end{align}
where we used (\ref{eq:lambda0}) and (\ref{eq:lambda2}). Merging the preceding inequality with (\ref{eq:EcBc}) and (\ref{eq:fLowerEcB}) we get for $n\ge n_0$
\begin{align}
	\int_{(E_{\varepsilon}(\ddot{f}(\mathbf{0})))^c}e^{-nf(\mathbf{t})}d\mathbf{t}&=
	\int_{(E_{\varepsilon}(\ddot{f}(\mathbf{0})))^c}e^{-(n-n_1)f(\mathbf{t})}e^{-n_1f(\mathbf{t})}d\mathbf{t}\notag\\
		&\le\notag\frac74e^{\xi n_1/2}n^{-({d}/2+\alpha)}\int_{(E_{\varepsilon}(\ddot{f}(\mathbf{0})))^c}e^{-n_1f(\mathbf{t})}d\mathbf{t}\notag\\
		&\le\notag\frac74e^{\xi n_1/2}n^{-({d}/2+\alpha)}I(n_1),\notag
\end{align}
which completes the proof.\\
$\square$

\begin{remark}
	Clearly, the threshold $7/4$ in the condition (\ref{eq:lambda0}) ensues from the range of validity of the inequality $e^x\le 1+x+x^2$ and this arbitrary value may impose a strong restriction on the value $n_0$ and thus may limit the range of applicability of the result. When the condition (\ref{eq:lambda0}) is more restrictive than (\ref{eq:lambda2}) one may apply instead an inequality of a form $e^x\le1+x+(1+a)x^2$ valid for $x\le x_a$ with $a>0$ (and hence $x_a>7/4$) chosen such that (\ref{eq:lambda0}) holds for all $n\ge n_2$. This will, however, be done at a cost of enlarging the constants $K_1$, $K_{\alpha,2}$ and $K_u$ which should then be respectively rewritten as $aK_1$, $aK_{\alpha,2}$ and $4/7x_a$. Such a trade-off between the quality of the upper bound in (\ref{Lap_appr_statement_Upper}) and its range of applicability is illustrated in example \ref{example_d=2}.
	On the other hand a completely opposite situation may arise, where the above-mentioned constants in (\ref{Lap_appr_statement_Upper}) are unsatisfactory but the left hand side of (\ref{eq:lambda0}) leaves plenty of space for adjustment. In such a situation one may take a suitable $a\in(-1/2,0)$ and improve $K_1$, $K_{\alpha,2}$ and $K_u$ at a cost (if any) of enlarging $n_0$.
\end{remark}

The form of the upper bound in theorem \ref{th:muldimlapl} involves two possible leading terms of the main error. Depending on $\alpha$ the rate of the main error term is either $n^{-\alpha/2}$ or $n^{-1}$. Inglot and Majerski proved that the rate $n^{-\alpha/2}$ is not underestimated when third order partial derivatives all vanish at $\mathbf{0}$, and hence the $n^{-1}$ term would be absent in our theorem \ref{th:muldimlapl}. The following example shows that it is also the case without this condition for $\alpha<2$.

\begin{example}
	Let $f(\boldsymbol t)=\sum_{i=1}^dt_i^2+\sum_{i=1}^dt_i^3+\sum_{i=1}^d|t_i|^{3+\gamma}$, where $\gamma\in(0,1)$. The assumtions of theorem \ref{th:muldimlapl} all hold. In particular for every $\mathbf{t}\in\mathbb{R}^d$ 
\[
\left|f(\mathbf{t})-\frac{1}{2}d^2f(\boldsymbol0,\boldsymbol t)-\frac16d^3f(\mathbf{0},\mathbf{t})\right|\le ||\mathbf{t}||^{3+\gamma},
\]
which means that (A1) holds for each $r>0$ with $\alpha=\gamma+1$ and $C=1$. The condition (A2) holds with $\delta=r$. Moreover, $\ddot{f}(\mathbf{0})={\rm diag}(2,\ldots,2)'$, so $\lambda_{\rm min}=2$. All the unmixed third order partial derivatives at $\mathbf{0}$ are equal 6, mixed being equal $0$, which yields $d^3f(\boldsymbol0,\boldsymbol t)/6=\sum_{i=1}^dt_i^3$.

We have
\begin{equation}\label{eq:ex1}
	\int_{\mathbb{R}^d}e^{-n f(\mathbf{t})}d\mathbf{t}=\left(\int_{-\infty}^{\infty}e^{-n(x^2+x^3+|x|^{3+\gamma})}dx\right)^d=(i(n))^d.
\end{equation}
Using the inequalities $e^{-x}\ge1-x$ $(x\in\mathbb{R})$ and $e^{-x}\le1-x+x^2/2$ ($x>0$) we get respectively
\[
	i(n)\ge\int_{-\infty}^{\infty}e^{-nx^2}\left(1-nx^3-n|x|^{3+\gamma})\right)dx=\sqrt{\frac\pi n}\left(1-\frac{\Gamma(2+\frac{\gamma}{2})}{\sqrt{\pi}}\cdot n^{-(1+\gamma)/2}\right),
\]
and  
\begin{align}
	\int_0^\infty e^{-n(x^2+x^3+x^{3+\gamma})}dx &\le 
\frac12\sqrt{\frac{\pi}{n}}\left(1-\frac{\Gamma(2+\frac{\gamma}{2})}{\sqrt{\pi}}\cdot n^{-(1+\gamma)/2}-\frac{1}{\sqrt{\pi n}}\right.\notag\\
		&+\left.\frac{15}8{n}^{-1}+\frac{2\Gamma((7+\gamma)/2)}{\sqrt{\pi}}n^{-(1+\gamma/2)}+
\frac{\Gamma(7/2+\gamma)}{\sqrt{\pi}}n^{-(1+\gamma)}\right)\notag.
\end{align}
Put $a_n:=\sqrt{2\log n/n}$. We have
\begin{align}
	\int_{-\infty}^{-a_n}e^{-n(x^2+x^3+(-x)^{3+\gamma})}dx&=\int_{a_n}^{\infty}e^{-(n-1)(x^2-x^3+x^{3+\gamma})}e^{-(x^2-x^3+x^{3+\gamma})}dx\notag\\
					& \le e^{-na_n^2}e^{na_n^3-na_n^{3+\gamma}}\int_0^{\infty}e^{-(x^2-x^3+x^{3+\gamma})}dx\notag\\
					& \le C_1n^{-2},\notag
\end{align}
with $C_1=\exp(0.3)\int_0^{\infty}e^{-(x^2-x^3+x^{3+\gamma})}dx$. Moreover, since $x^3-x^{3+\gamma}$ is increasing for sufficiently small $x>0$, it holds that for $n$ large enough and $x\le a_n$
\[
	n(x^3-x^{3+\gamma})\le\frac{(2\log n)^{3/2}}{n^{1/2}}\left(1-\left(\frac{2\log n}{n}\right)^{\gamma/2}\right)<7/4.
\]
Thus, using the inequality $e^{x}\le1+x+x^2$, $(x\le7/4)$,
\begin{align}
	\int_{0}^{a_n}e^{-nx^2}e^{nx^3-nx^{3+\gamma}}dx&\le \int_{0}^{a_n}e^{-nx^2}\left(1+nx^3-nx^{3+\gamma}+n^2(x^3-x^{3+\gamma})^2\right)dx\notag\\
					&\le \int_0^{\infty}e^{-nx^2}(1+nx^3+n^2x^6)dx-n\int_{0}^{a_n}e^{-nx^2}x^{3+\gamma}dx\notag\\
					&=\notag \frac12\sqrt{\frac{\pi}{n}}\left(1+\frac{1}{\sqrt{\pi n}}+\frac{15}8n^{-1}-\frac{\Gamma(2+\frac{\gamma}{2})-\Gamma(2+\frac{\gamma}{2},2\log n)}{\sqrt{\pi}}\cdot n^{-(1+\gamma)/2}\right),\notag
\end{align}
where $\Gamma(\cdot,\cdot)$ denotes the incomplete gamma function. In view of a fact $\Gamma(k,2\log n)\sim (2\log n)^{k-1}n^{-2}$ as $n\to \infty$ (\cite{GR}), we infer that with $C_2<15/8$ for every large $n$ it holds
\[
	i(n)\le\sqrt{\frac{\pi}{n}}\left(1-\frac{\Gamma(2+\frac{\gamma}{2})}{\sqrt{\pi}}\cdot n^{-(1+\gamma)/2}+C_2n^{-1}\right).
\]

Hence from (\ref{eq:ex1}) and the multinomial theorem we conclude that
\begin{align}
	I(n)= \left(\frac{\pi}{n}\right)^{d/2}\left(1-\frac{d\Gamma(2+\frac{\gamma}{2})}{\sqrt{\pi}}\cdot n^{-(1+\gamma)/2}+O(n^{-1})\right).
\end{align}
\end{example}
\vspace{.5cm}

\begin{theorem}\label{th:muldimlapl_g}
    Let $\Omega$ and $f$ satisfy assumptions of theorem \ref{th:muldimlapl}. Let $g:\Omega\to\mathbb{R}$ be a measurable function and $g(\boldsymbol0)\neq0$. Suppose $g\in {\rm C}^{1,1}(B_r)$, so that\\
		
        {\rm (A3)}\quad there exists $M>0$ such that for every $\boldsymbol t\in B_r$ 
		$$\left|g(\boldsymbol t)-g(\boldsymbol 0)-\dot{g}(\boldsymbol 0)'\boldsymbol t\right|\le M\left\|t\right\|^2.$$
		
Let $J(n)=\int_{\Omega}e^{-n f(\boldsymbol t)}g(\boldsymbol t)dt$ and suppose there exists $n_3>0$ such that $J(n)$ is finite for all $n\geqslant n_3$. Then for every $n\geqslant n_4$, where $n_4$ is given by (\ref{eq:lambda4}), it holds
    \[
    J(n)=\frac{g(\boldsymbol0)}{\sqrt{\mbox{\rm det} \ddot{f}(\boldsymbol0)}}\left(\frac{2\pi}{n}\right)^{d/2}\left(1+E(n)\right),
    \]
    where
    \begin{align}
     E(n)&\geqslant-\frac{K_{\alpha,1}}{n^{\alpha/2}}-\frac{K_2+K_3}{n}-\frac{K_{\alpha,3}}{n^{1+\alpha/2}}-\frac{K_{ul}}{n^\alpha}-\frac{K_l}{n^{1+\alpha}};
    \notag\\
	E(n)&\leqslant \frac{K_{\alpha,1}}{n^{\alpha/2}}+\frac{K_1+K_2+K_3}{n}+\frac{K_{\alpha,2}+K_{ul}}{n^{\alpha}}+
        \frac{K_{\alpha,3}}{n^{1+\alpha/2}}+\frac{K_4}{n^2}+\frac{K_{\alpha,5}}{n^{1+\alpha}}+\frac{K_{\alpha,6}}{n^{(3+\alpha)/2}},\notag
    \end{align}
    with $K_{\alpha,1},\;K_{\alpha,2},\;K_1$ and $K_l$ are given by (\ref{constant_K}), (\ref{constant_K_1}) and (\ref{constant_K_l}), respectively, and
	\[
		K_2=\frac{M\,(2/\lambda_{\rm{min}})(d/2)_1}{\left|g(\boldsymbol0)\right|},\;
		K_3=\frac{D\|\dot{g}(\boldsymbol0)\|\,(2/\lambda_{\rm{min}})^2(d/2)_2}{|g(\boldsymbol0)|},
	\]
	\[        
		K_{\alpha,3}=\frac{CM\,(2/\lambda_{\rm{min}})^{2+\alpha/2}(d/2)_{2+\alpha/2}}{\left|g(\boldsymbol0)\right|},
        \;K_4=\frac{MD^2(2/\lambda_{\rm{min}})^4(d/2)_4}{|g(\boldsymbol0)|},
    \]
	\[
        K_{\alpha,5}=\frac{C^2M\,(2/\lambda_{\rm{min}})^{1+\alpha}(d/2)_{1+\alpha}}{\left|g(\boldsymbol0)\right|},
        \;K_{\alpha,6}=\frac{2CD\|\dot{g}(\boldsymbol0)\|\,(2/\lambda_{\rm{min}})^{(3+\alpha)/2}(d/2)_{(3+\alpha)/2}}{\left|g(\boldsymbol0)\right|},
    \]
	\[
        K_{ul}=\frac{7\sqrt{\mbox{\rm det}\,\ddot{f}(\boldsymbol0)}}{4\left|g(\boldsymbol0)\right|(2\pi)^{d/2}}e^{n_3\xi/2}\int_{\mathbb{R}^d}e^{-n_3 f(\boldsymbol t)}|g(\boldsymbol t)|d\boldsymbol t.
    \]
    Moreover, $n_4\ge n_0$ is the smallest number such that the following inequality holds true
\begin{equation}\label{eq:lambda4}
	\frac{M(d+2\alpha)}{\lambda_{\rm min}}\frac{\log n}{n}+\|\dot{g}(\boldsymbol0)\|\sqrt{\frac{d(d+2\alpha)}{\lambda_{\rm min}}}\cdot\sqrt{\frac{\log n}{n}}\le g(\boldsymbol0).
\end{equation}
\end{theorem}

\noindent{\bf Proof.} The proof is a slight extension that of theorem 2.3 in \cite{IM} and an extension of our equality (\ref{eq:int_d3}). The condition (\ref{eq:lambda4}) insures that for every $t\in B_r$ the both bounds for $g$ in the inequalities
\[
	g(\boldsymbol 0)+\dot{g}(\boldsymbol 0)'\boldsymbol t-M\|t\|^2 \le g(\boldsymbol t)\le g(\boldsymbol 0)+\dot{g}(\boldsymbol 0)'\boldsymbol t+M\|t\|^2
\]
are nonnegative. 
$\square$\\

A consequence of the above theorem is the following corollary, which provides the rate $n^{-1}$ under the weakest regularity assumption on $f$ and $g$.
 
\begin{corollary}\label{cor}
	Let the conditions (L1)-(L4) all hold. Let moreover $f\in{\rm C}^{2,1}(\boldsymbol t_0)$ and $g\in{\rm C}^{1,1}(\boldsymbol t_0)$. Then it holds
	\[
		J(n)=\frac{e^{-n f(\boldsymbol t_0)}g(\boldsymbol t_0)}{\sqrt{\mbox{\rm det} \ddot{f}(\boldsymbol t_0)}}\left(\frac{2\pi}{n}\right)^{d/2}\left(1+O(n^{-1})\right).
	\]
\end{corollary}
\vspace{.5cm}

\section{Application}

\begin{example}
Following \cite{deB} consider the sum
\[
	S(s,n)=\sum_{k=0}^{2n}(-1)^{k+n}{2n \choose k}^s,
\]
where $n$ and $s$ are two positive integers. This formula is known in an explicit form for at most $s=3$ (called Dixon's identity) and no such extension exists for $s\ge4$. Here we consider an approximation of $S(s,n)$ as $n\to\infty$ for $s$ fixed. To this end we shall use its following integral representation \cite[p. 73]{deB}
\[
	S(d+1,n)=\frac{2^{2n(d+1)}}{\pi^d}\cdot\int_{(-\pi/2,\pi/2)^d}\left\{\sin\!\left(\sum_{i=1}^d \varphi_i\right)\prod_{i=1}^d\cos\varphi_i\right\}^{2n} d\boldsymbol\varphi,
\]
with $\boldsymbol\varphi=(\varphi_1,\ldots,\varphi_d)'$. As the integrand is symmetrical over values of $\boldsymbol\varphi$ corresponding to $\sum_{i=1}^d \varphi_i>0$ or $<0$, in order to investigate the asymptotic behavior of $S(d+1,n)$ it suffices to consider the integral 
\[
	T(d,n):=\int_{\Omega'} e^{-2n h(\boldsymbol\varphi)}d\boldsymbol\varphi,
\] 
where $\Omega'=\{\boldsymbol\varphi\in\mathbb{R}^d:~\varphi_i<\pi/2,~\varphi_1+\ldots\varphi_d>0\}$, and $h(\boldsymbol\varphi)=-\log\sin(\sum_{i=1}^d \varphi_i)-\sum_{i=1}^d\log\cos\varphi_i$. The unique minimizer of $h$ in $\Omega'$ is $\hat{\boldsymbol\varphi}=(\hat\alpha,\ldots,\hat\alpha)'$ with $\hat\alpha=\pi/2(d+1)$ and $h(\hat{\boldsymbol\varphi})=-(d+1)\log\cos(\hat\alpha)$, for $\sin(d\hat\alpha)=\cos(\hat\alpha)$. In order to apply theorem \ref{th:muldimlapl} write
\[
	T(d,n)=\cos^{2n(d+1)}(\hat{\alpha})\int_\Omega e^{-2n f(\boldsymbol\varphi)}d\boldsymbol\varphi,
\]
where 
\[
	\Omega=\{\boldsymbol\varphi\in\mathbb{R}^d:~\varphi_i<d\hat{\alpha},~\varphi_1+\ldots\varphi_d>-d\hat{\alpha}\}
\] 
and $f(\boldsymbol\varphi)=h(\boldsymbol\varphi+\hat{\boldsymbol\varphi})-h(\hat{\boldsymbol\varphi})=
\log\cos(\sum_{i=1}^d\varphi_i-\hat\alpha)-
\sum_{i=1}^d\log\cos(\varphi_i+\hat\alpha)-
h(\hat{\boldsymbol\varphi})$. Now $\mathbf{0}$ minimizes $f$ over $\Omega$ and 
$\dot{f}(\mathbf{0})=0$. We also have $\ddot{f}(\mathbf{0})=[(\delta_{ij}+1)\cos^{-2}(\hat{\alpha})]_{i,j=1,\ldots,d}$, where $\delta_{ij}$ denotes the Kronecker delta, hence $\det(\ddot{f}(\mathbf{0}))=(d+1)\cos^{-2d}(\hat{\alpha})$ and $\lambda_{\rm{min}}=\cos^{-2}(\hat{\alpha})$. Moreover, 
\[
	\frac{\partial ^3f(\boldsymbol0)}{\partial t_i\partial t_j\partial t_k}=-2(1-\delta_{ij}\delta_{jk})(\cos^{-3}\cdot\sin)(\hat{\alpha}).
\]

For fixed $\eta\in(0,1)$ let us define the ball $B_{\eta\sqrt{d}\hat{\alpha}}$. It is easy to verify that $B_{\eta\sqrt{d}\hat{\alpha}}\subset \Omega_{\eta}\subset \Omega$, where $\Omega_\eta=\{\boldsymbol\varphi\in\mathbb{R}^d:~|\varphi_i|\le\eta\sqrt{d}\hat{\alpha},~|\varphi_1+\ldots\varphi_d|\le\eta d\hat{\alpha}\}$. An application of Taylor's theorem yields for $\boldsymbol\varphi\in B_{\eta\sqrt{d}\hat{\alpha}}$
\[
	\left|f(\boldsymbol\varphi)-\frac12d^2f(\boldsymbol0,\boldsymbol \varphi)-\frac16d^3f(\boldsymbol0,\boldsymbol\varphi)\right|=
		\frac1{24}\left|d^4(\boldsymbol\xi,\boldsymbol\varphi)\right|,
\]
where $\boldsymbol\xi\in B_{\eta\sqrt{d}\hat{\alpha}}\subset \Omega_\eta$. We have
\[
	\frac{\partial ^4f(\boldsymbol\xi)}{\partial t_i\partial t_j\partial t_k\partial t_l}=2\left(u\left(\sum_{i=1}^d\xi_i-\hat{\alpha}\right)
		+\delta_{ij}\delta_{jk}\delta_{kl}\cdot u\left(\xi_l+\hat{\alpha}\right)\right),
\]
where $u:=(\cos^{-4}\cdot(1+2\sin^2))$. It is seen that $u(\sum_{i=1}^d\varphi_i-\hat\alpha)\le u((\eta d+1)\hat\alpha)$ and $u(\varphi_i+\hat\alpha)\le u((\eta\sqrt{d}+1)\hat\alpha)$ on $\Omega_\eta$. One thus obtains for every $\boldsymbol\varphi\in B_{\eta\sqrt{d}\hat{\alpha}}$
\begin{equation}\label{eq:Ex3C}
	\left|f(\boldsymbol\varphi)-\frac12d^2f(\boldsymbol0,\boldsymbol \varphi)-\frac16d^3f(\boldsymbol0,\boldsymbol\varphi)\right|\le
	C\|\boldsymbol\varphi\|^4,
\end{equation}
where 
\[
	C=\frac{d^2}{12}\left\{
				\frac{1+2\sin^2}{\cos^4}\left(\frac{(\eta d+1)\pi}{2(d+1)}\right)+\frac{1+2\sin^2}{\cos^4}\left(\frac{(\eta\sqrt{d}+1)\pi}{2(d+1)}\right)
					\right\},
\]
thus, the condition (A1) holds. To verify (A2) note that $f$ is convex on $\Omega$ since $\ddot{f}(\boldsymbol\varphi)=(\cos^{-2}(\sum_{i=1}^d\varphi_i-\hat{\alpha})+\delta_{ij}\cos^{-2}(\varphi_j+\hat{\alpha}))_{i,j=1,\ldots,d}$ is postive definite as a matrix of the quadratic form
\[
	\cos^{-2}(\varphi_1+\hat{\alpha})x_1^2+\ldots+\cos^{-2}(\varphi_d+\hat{\alpha})x_d^2+\cos^{-2}(\sum_{i=1}^d\varphi_i-\hat{\alpha})(x_1+\ldots +x_d)^2.
\]
We thus may take $\delta=r=\eta\sqrt{d}\hat{\alpha}$. Finding a reasonable estimation for $\Delta$ in such a general setting seems a difficult task. One way to do this is to use the (A1) condition, from which $f(\boldsymbol\varphi)\ge\boldsymbol\varphi\ddot{f}(\boldsymbol 0)\boldsymbol\varphi'/2+d^3f(\boldsymbol 0,\boldsymbol\varphi)/6-C\|\boldsymbol\varphi\|^4$ on $B_r$ and because $|d^3f(\boldsymbol 0,\boldsymbol\varphi)/6|\le D\|\boldsymbol\varphi\|^3$ and $\boldsymbol\varphi\ddot{f}(\boldsymbol 0)\boldsymbol\varphi'\ge\lambda_{\rm min}\|\boldsymbol\varphi\|^2$ we infer that for every $t\in B_r$
\[
	f(\boldsymbol\varphi)\ge\frac{\lambda_{\rm min}}2r^2-Dr^3-Cr^4.
\]
Hence, in view of convexity of $f$ one may take $\Delta:=\lambda_{\rm min}r^2/2-Dr^3-Cr^4$, provided $\Delta>0$, i.e. $\eta$ is small enough.

Denoting $I(n):=\int_\Omega e^{-2nf(\boldsymbol\varphi)}d\boldsymbol\varphi$, we finally conclude that for every $\eta\in(0,1)$ such small that $\Delta$ defined above is positive, with $r:=\eta\sqrt{d}\hat\alpha$ and for every large enough $n\ge1$ there hold
\begin{align}
	I(n)&\ge\frac{\cos^d(\hat\alpha)}{\sqrt{d+1}}\left(\frac{\pi}n\right)^{d/2}
					\left\{1-\frac{K_{\alpha,1}}{2n}-\frac{K_l}{8n^3}\right\},\notag\\
	I(n)&\le\frac{\cos^d(\hat\alpha)}{\sqrt{d+1}}\left(\frac{\pi}n\right)^{d/2}
					\left\{1+\frac{K_{\alpha,1}+K_1}{2n}+\frac{K_{\alpha,2}+K_u}{4n^2}\right\},\notag
\end{align}
with $K_{\alpha,1}=C\cos^4(\hat\alpha)d(d+2)$, $K_1=\sin^2(\hat\alpha)d^4(d+2)(d+4)$, 
$K_{\alpha,2}=C^2\cos^8(\hat\alpha)d(d+2)(d+4)(d+6)$, $K_l$ given by (\ref{constant_K_l}) and 
$K_u=(7/4)\sqrt{d+1}\cos^{-d}(\hat\alpha)(2\pi)^{-d/2}I(1)e^{\xi/2}$, where $I(1)=(2^d-1)(\pi/4)^d$.

Recalling that $\hat\alpha=\pi/2(d+1)$ one sees consistency in the behavior of our main error term of the upper and lower bound as $d\to\infty$, for $K_{\alpha,1}\asymp d^4$ and $K_1\asymp d^4$. 
\end{example}

\begin{example}[continued]\label{example_d=2}
	We shall consider the case $s=3$ and compare our estimations with those obtained previously in \cite{MCW}. For $d=2$ $\Delta>0$ whenever $\eta\le0.36$. Taking example value $\eta=1/3$ we find $C=7.7$, $K_{\alpha,1}/2=17.4$ and $K_1/2=64$. Hence for every large $n$ we have
	\[
		S(3,n)=\frac{3^{3n+\frac12}}{2\pi n}\left(1+\frac{A}{n}+O\!\left(\frac1{n^2}\right)\right),
	\]
	where $A\le 81.4$. Note that taking $\eta$ small enough one can near the constant $C$ arbitrarily close to $C^{(best)}=16/9$, and consequently reach $A^{(best)}=68$. As we shall see this is a very poor estimate for $A$, which is a consequence of a poor estimation of $C$.

We shall now adopt the method used by McClure and Wong (\cite{MCW}), who proposed prior elimination of cross product term $xy$ in the Taylor expansion of $f(\varphi_1,\varphi_2)$ via a linear transformation. This can be achieved using the Cholesky decomposition of $\ddot{f}(0,0)$. Thus, following \cite[p. 388]{MCW} and applying to our $I(n)$ the change of variables
\[
	\varphi_1=(\sqrt{3}/2)x-(1/2)y,~~~~\varphi_2=y
\]
one obtains
\[
	S(3,n)=\frac{{3}^{3n+\frac12}}{\pi^2}\int_{\Omega}e^{-2nf(x,y)}dxdy,
\]
where
\begin{align}
	f(x,y)= &-\log\cos \left(\frac{\sqrt3}{2}x-\frac12y+\frac{\pi}6\right)-\log\cos \left( -\frac{\sqrt3}2x-\frac12y+\frac{\pi}6 \right)\notag\\
			&-\log\cos \left( y+\frac\pi6 \right) +3\log\frac{\sqrt3}2,\notag
\end{align}
and $\Omega=\left\{(x,y)\in\mathbb{R}^2:\frac{\sqrt3}2x-\frac12y<\frac\pi3,~ y<\frac\pi3,~ \frac{\sqrt3}2x+\frac12y>-\frac\pi3\right\}$.

Take $r\in(0,\pi/3)$ and note that $B_r\subset\mathbb{R}^2$ is enclosed by the hexagonal set $\Omega_r:=\{(x,y)\in\mathbb{R}^2:~|\frac{\sqrt3}2x-\frac12y|\le r,~ |y|\le r,~ |\frac{\sqrt3}2x+\frac12y|\le r\}$. Considering the univariate function $-\log\cos(t+\pi/6)$ whose Maclaurin expansion is of a form
\[
	-\log\cos\left(t+\frac{\pi}6\right)=\frac1{\sqrt3}t+\frac23t^2+\frac4{9\sqrt3}t^3+\frac29t^4+\ldots
\]
we can apply the Taylor theorem with remainder to each of the three first terms defining $f(x,y)$ and obtain on $B_r$
\begin{align}
	f(x,y)-x^2&-y^2-\frac{\sqrt{3}}9(y^3-3x^2y)\notag\\
				  &=\frac1{12}\left\{ \left(\frac{\sqrt3}{2}x-\frac12y\right)^4 u(\xi_1)+\left(-\frac{\sqrt3}{2}x-\frac12y\right)^4 u(\xi_2)+y^4 u(\xi_3)\right\},\notag
\end{align}
where $u(t):=(\cos^{-4}(1+2\sin^2))	(t+\pi/6)$ and $\xi_i\in(-r,r)$, for $i=1,2,3$. Thus, using the relation
\[
	(({\sqrt3}/{2})x-(1/2)y)^4+(({\sqrt3}/{2})x+(1/2)y)^4+y^4=(9/8)(x^2+y^2)^2
\]
we infer 
\begin{equation}\label{ex:Cbound}
	\left|f(x,y)-x^2-y^2-\frac{\sqrt{3}}9(y^3-3x^2y)\right|\le C(x^2+y^2)^2,
\end{equation}
with $C=(3/32)(\cos^{-4}(1+2\sin^2))(r+\pi/6)$. Following \cite{MCW} put $r^2=\pi^2/108$ and note that the condition (A1) holds with $C=0.9238$, $\alpha=2$ and $\lambda_{\rm{min}}=2$. Moreover, $\delta=r$, $f$ is convex and reasoning as in the previous example we find $\Delta=0.06863$. Hence we have
\begin{align}
	S(3,n)&\ge\frac{{3}^{3n+1/2}}{2\pi n}
					\left\{1-\frac{K_{\alpha,1}/2}{n}-\frac{K_l/8}{n^3}\right\},\label{eq:S3nApprLo}\\
	S(3,n)&\le\frac{{3}^{3n+1/2}}{2\pi n}
					\left\{1+\frac{(K_{\alpha,1}+K_1)/2}{n}+\frac{(K_{\alpha,2}+K_u)/4}{n^2}\right\},\label{eq:S3nApprUp}
\end{align}
where $K_{\alpha,1}/2=0.9238$, $(K_{\alpha,1}+K_1)/2=1.072$, $K_l/8=0.1355$, $K_{\alpha,2}=20.48$ and $(K_{\alpha,2}+K_u)/4=5.439$.

Using their device McClure and Wong \cite[eq. (5.37)-(5.39)]{MCW} proved that (in our notation and emending an obvious slip of the pen in the leading term in (5.37); cf. \cite[p.75]{deB})
\[
	S(3,n)=\frac{{3}^{3n+1/2}}{2\pi n}
					\left\{1+\tilde{E}_1(n)+\tilde{E}_2(n)\right\},
\]
where $|\tilde{E}_1(n)|\le 1.8245/n$ and $|\tilde{E}_2(n)|\le (7/3)e^{-n\pi^2/72}$. We conclude that our main error term $1.072$ significantly improves the known estimation from \cite{MCW}. We also note that for every $n\ge2$ the relative error given in our bounds is better estimated than the compared one. For instance, when $n=1$ their result implies that the error is within $\mp 3.9$ whereas our gives the interval $(-1.1,6.5)$. For $n=2$ these are $\mp2.7$ and $(-0.48,1.9)$; for $n=5$: $\mp1.5$ and $(-0.18,0.43)$; for $n=10$: $\mp0.77$ and $(-0.093,0.16)$ while for $n=100$: $\mp0.018$ and $(-0.0092,0.011)$.

It is however well known that the approach used in this paper fails to indicate a satisfactory value for $n_0$, i.e. the minimum value of $n$ insuring the validity of inequalities in (\ref{eq:S3nApprLo}) and (\ref{eq:S3nApprUp}) (cf. \cite[p. 11]{IM}). While the right hand side of (\ref{eq:lambda2}) already indicates $n_0\ge240$, the requirement (\ref{eq:lambda0}) yields $n_0=1479$. Meanwhile, numerical simulations not presented here showed that the inequalities work for every positive integer $n$. In order to make the condition (\ref{eq:lambda0}) less restrictive than (\ref{eq:lambda2}) one could use a weaker inequality $e^x\le1+x+2.2x^2$ $(x<3.39)$ in the proof of the theorem \ref{th:muldimlapl} yielding $n_0=240$ at a cost of enlarging the constants $(K_{\alpha,1}+K_1)/2$ and $(K_{\alpha,2}+K_u)/4$ to $1.2497$ and $11.5833$, respectively. On the other hand, even if one could find a better estimation for $\Delta$, for a given $r^2=\pi^2/108$ the restriction imposed by (\ref{eq:lambda2}) gives at least $n_0\ge169$, for $\xi\ge r^2\lambda_{\rm min}$.
\end{example}




\end{document}